\input amstex
\documentstyle{amsppt}
\magnification=\magstephalf \hsize = 6.5 truein \vsize = 9 truein
\vskip 3.5 in

\rightheadtext\nofrills {ALGORITHM THAT DECIDES TRANSLATION
EQUIVALENCE} \NoBlackBoxes \TagsAsMath

\def\label#1{\par%
        \hangafter 1%
        \hangindent .75 in%
        \noindent%
        \hbox to .75 in{#1\hfill}%
        \ignorespaces%
        }

\newskip\sectionskipamount
\sectionskipamount = 24pt plus 8pt minus 8pt
\def\sectionskip{\vskip\sectionskipamount}
\define\sectionbreak{%
        \par  \ifdim\lastskip<\sectionskipamount
        \removelastskip  \penalty-2000  \sectionskip  \fi}
\define\section#1{%
        \sectionbreak   
        \subheading{#1}%
        \bigskip
        }

\redefine\qed{{\unskip\nobreak\hfil\penalty50\hskip2em\vadjust{}\nobreak\hfil
    $\square$\parfillskip=0pt\finalhyphendemerits=0\par}}

        
        \let    \< = \langle
        \let    \> = \rangle

\define\op#1{\operatorname{\fam=0\tenrm{#1}}} 
        \define         \a              {\alpha}
        \redefine       \b              {\beta}
        \redefine       \d              {\delta}
        \redefine       \D              {\Delta}
        \define         \e              {\varepsilon}
        \define         \E              {\op {E}}
        \define         \g              {\gamma}
        \define         \G              {\Gamma}
        \redefine       \l              {\lambda}
        \redefine       \L              {\Lambda}
        \define         \n              {\nabla}
        \redefine       \var            {\varphi}
        \define         \s              {\sigma}
        \redefine       \Sig            {\Sigma}
        \redefine       \t              {\tau}
        \define         \th             {\theta}
        \redefine       \O              {\Omega}
        \redefine       \o              {\omega}
        \define         \z              {\zeta}
        \define         \k              {\kappa}
        \redefine       \i              {\infty}
        \define         \p              {\partial}
        \define         \vsfg           {\midspace{0.1 truein}}

\topmatter
\title An algorithm that decides translation equivalence in a free
group of rank two
\endtitle
\author Donghi Lee
\endauthor

\address {Department of Mathematics, Pusan National University, Jangjeon-Dong, Geumjung-Gu, Pusan 609-735, Korea}
\endaddress

\email {donghi\@pusan.ac.kr}
\endemail


\abstract {Let $F_2$ be a free group of rank $2$. We prove that
there is an algorithm that decides whether or not, for given two
elements $u, v$ of $F_2$, $u$ and $v$ are translation equivalent
in $F_2$, that is, whether or not $u$ and $v$ have the property
that the cyclic length of $\phi(u)$ equals the cyclic length of
$\phi(v)$ for every automorphism $\phi$ of $F_2$. This gives an
affirmative solution to problem F38a in the online version
(http://www.grouptheory.info) of [1] for the case of $F_2$.}
\endabstract
\endtopmatter


\document
\baselineskip=24pt

\heading 1. Introduction
\endheading

Let $F_n$ be the free group of rank $n \ge 2$ on the set $\Sigma$.
As usual, for a word $v$ in $F_n$, $|v|$
denotes the length of the reduced word over $\Sigma$ representing
$v$. A word $v$ is called {\it cyclically reduced} if all its
cyclic permutations are reduced. A {\it cyclic word} is defined to
be the set of all cyclic permutations of a cyclically reduced
word. By $[v]$ we denote the cyclic word associated with a word
$v$. Also by $\|v\|$ we mean the length of the cyclic word $[v]$
associated with $v$, that is, the number of cyclic permutations of
a cyclically reduced word which is conjugate to $v$. The length
$\|v\|$ is called the {\it cyclic length} of ~$v$.

In [2], Kapovich-Levitt-Schupp-Shpilrain introduced and studied in
detail the notion of translation equivalence in free groups. The
following definition is a combinatorial version of translation
equivalence:

\proclaim{Definition 1.1 {\rm [2, Corollary 1.4]}} Two words $u, v
\in F_n$ are called {\it translation equivalent in $F_n$} if the
cyclic length of $\phi(u)$ equals the cyclic length of $\phi(v)$
for every automorphism $\phi$ of $F_n$.
\endproclaim

Several different sources of translation equivalence in free
groups were provided by Kapovich-Levitt-Schupp-Shpilrain [2] and
Lee [3]. Pointing out in [2] that hyperbolic equivalence in
surface groups (cf. [5]) and character equivalence in free groups
are algorithmically decidable, Kapovich-Levitt-Schupp-Shpilrain
raised the question about the existence of an algorithm which
decides translation equivalence in free groups.

The purpose of the present paper is to prove that translation
equivalence is algorithmically decidable in $F_2$.

\proclaim {Theorem 1.2} There exists an algorithm that decides
whether or not, for given two elements $u, \, v$ of $F_2$, $u$ and
$v$ are translation equivalent in $F_2$.
\endproclaim

In conclusion as will be shown in Section 3,
the algorithm in the statement of Theorem 1.2 is as follows.

\proclaim {Algorithm} Let $F_2=\<x,y\>$, and let $\O$ be the set
of all chains of Whitehead automorphisms of $F_2$ of the form
either
$$(\{y\}, x)^{m_k}(\{x\}, y)^{l_k} \cdots (\{y\},
x)^{m_1}(\{x\}, y)^{l_1}$$ or
$$(\{y\}, x^{-1})^{m_k}(\{x\}, y^{-1})^{l_k} \cdots (\{y\},
x^{-1})^{m_1}(\{x\}, y^{-1})^{l_1},$$ where $k \in \Bbb{N}$, each
$l_i, \, m_i \ge 0$ and $\sum_{i=1}^k (l_i+m_i) \le 2\|u\|+3$.
Then $\O$ is clearly a finite set. Check if
$\|\psi(u)\|=\|\psi(v)\|$ for every $\psi \in \O$. If so, conclude
that $u$ and $v$ are translation equivalent in $F_2$; otherwise
conclude that $u$ and $v$ are not translation equivalent in $F_2$.
\endproclaim

Here, as in [4], a {\it Whitehead automorphism} $\s$ of $F_n$ is
defined to be an automorphism of one of the following two types
(cf. [6]):

\roster \item"(W1)" $\s$ permutes elements in $\Sig^{\pm 1}$.
\item"(W2)" $\s$ is defined by a set $S \subset \Sig^{\pm 1}$ and
a letter $a \in \Sig^{\pm 1}$ with both $a, \, a^{-1} \notin S$ in
such a way that if $c \in \Sig^{\pm 1}$ then (a) $\s(c)=ca$
provided $c \in S$ and $c^{-1} \notin {S}$; (b) $\s(c)=a^{-1}ca$
provided both $c,\, c^{-1} \in S$; (c) $\s(c)=c$ provided both
$c,\, c^{-1} \notin S$.
\endroster
If $\s$ is of type (W2), we write $\s=(S, a)$. By $(\bar S,
a^{-1})$, we mean the Whitehead automorphism $(\Sig^{\pm 1} -S -
a^{\pm 1}, a^{-1})$. It is then easy to check that
$$(S, a)(w)=(\bar S,a^{-1})(w)
\tag 1.1
$$ for every cyclic word $w$ in $F_n$.

\heading 2. Preliminary Lemmas
\endheading

We begin this section by setting some notation. As in [2], if $w$
is a cyclic word in $F_n$ and $a, b \in \Sig^{\pm 1}$, we use
$n(w; a, b)$ to denote the total number of occurrences of the
subwords $ab$ and $b^{-1}a^{-1}$ in $w$. Then clearly $n(w;
a,b)=n(w; b^{-1}, a^{-1})$. Similarly we denote by $n(w;a)$ the
total number of occurrences of $a$ and $a^{-1}$ in $w$. Then again
clearly $n(w;a)=n(w; a^{-1})$. As in [4], for two automorphisms
$\phi$ and $\psi$ of $F_n$, by writing $\phi \equiv \psi$ we mean
the equality of $\phi$ and $\psi$ over all cyclic words in $F_n$,
that is, $\phi(w)=\psi(w)$ for every cyclic word $w$ in $F_n$.

From now on, let $F_2$ be the free group of rank $2$ on the set
$\{x,y\}$.

\proclaim {Lemma 2.1} Let $\a$ be a Whitehead automorphism of
$F_2$ of type (W2). Then exactly one of $\a \equiv 1$, $\a \equiv
(\{x\}, y)$, $\a \equiv (\{x\}, y^{-1})$, $\a \equiv (\{y\}, x)$
and $\a \equiv (\{y\}, x^{-1})$ is necessarily satisfied.
\endproclaim

\demo{Proof} Let $\a$ be a Whitehead automorphism of $F_2$ of type
(W2). By the definition of (W2), $\a$ is one of $(\{x\}, y)$,
$(\{x^{-1}\}, y)$, $(\{x^{\pm 1}\}, y)$, $(\{x\}, y^{-1})$,
$(\{x^{-1}\}, y^{-1})$, $(\{x^{\pm 1}\}, y^{-1})$, $(\{y\}, x)$,
$(\{y^{-1}\}, x)$, $(\{y^{\pm 1}\}, x)$, $(\{y\}, x^{-1})$,
$(\{y^{-1}\}, x^{-1})$ and $(\{y^{\pm 1}\}, x^{-1})$. Among these,
$(\{x^{\pm 1}\}, y)$, $(\{x^{\pm 1}\}, y^{-1})$, $(\{y^{\pm 1}\},
x)$ and $(\{y^{\pm 1}\}, x^{-1})$ play the same role as the
identity over every cyclic word in $F_2$. Moreover, by (1.1),
$(\{x^{-1}\}, y) \equiv (\{x\}, y^{-1})$, $(\{x^{-1}\}, y^{-1})
\equiv (\{x\}, y)$, $(\{y^{-1}\}, x) \equiv (\{y\}, x^{-1})$ and
$(\{y^{-1}\}, x^{-1}) \equiv (\{y\}, x)$ in $F_2$, thus proving
the lemma. \qed
\enddemo

Now for the rest of the paper, let $\s=(\{x\}, y)$ and $\t=(\{y\},
x)$ be Whitehead automorphisms of $F_2$. Then obviously
$\s^{-1}=(\{x\}, y^{-1})$ and $\t^{-1}=(\{y\}, x^{-1})$.

\proclaim {Lemma 2.2} In $F_2$, we have
$$
\alignat 3 \t^{-1}\s &\equiv \pi \t, &&\quad \t^{-1}\pi \equiv
\pi\s, &&\quad  \s^{-1}\pi \equiv \pi \t, \\
\s\t^{-1} &\equiv \pi \s^{-1}, &&\quad \t\pi \equiv \pi
\s^{-1}, &&\quad \s\pi \equiv \pi \t^{-1}, \\
\t\s^{-1} &\equiv \pi^{-1} \t^{-1}, &&\quad \t\pi^{-1} \equiv
\pi^{-1}\s^{-1}, &&\quad \s\pi^{-1} \equiv \pi^{-1} \t^{-1}, \\
\s^{-1}\t &\equiv \pi^{-1}\s, &&\quad \t^{-1}\pi^{-1} \equiv
\pi^{-1} \s, &&\quad \s^{-1}\pi^{-1} \equiv \pi^{-1} \t,
\endalignat
$$
where $\pi$ is a Whitehead automorphism of $F_2$ of type (W1) that
sends $x$ to $y$ and $y$ to $x^{-1}$.
\endproclaim

\demo{Proof} For the first equality, check that $(\pi
\t)^{-1}\t^{-1}\s=(\{y^{\pm 1}\}, x^{-1})(\{x^{\pm 1}\}, y) \equiv
1$ in $F_2$. In a similar way, the rest of the equalities can be
checked. \qed
\enddemo

\proclaim {Lemma 2.3} For every automorphism $\phi$ of $F_2$,
$\phi$ can be represented as $\phi \equiv \b \phi'$, where $\b$ is
a Whitehead automorphism of $F_2$ of type (W1) and $\phi'$ is a
chain of one of the forms
$$\aligned
(C1) \ & \phi' \equiv \t^{m_k}\s^{l_k} \cdots \t^{m_1}\s^{l_1} \\
(C2) \ & \phi' \equiv \t^{-m_k}\s^{-l_k} \cdots \t^{-m_1}\s^{-l_1}
\endaligned
$$
with $k \in \Bbb {N}$ and both $l_i, \, m_i \ge 0$ for every $i=1,
\dots, k$.
\endproclaim

\demo{Proof} By Whitehead's Theorem (cf. [6]) together with Lemma
2.1, an automorphism $\phi$ of $F_2$ can be expressed as
$$\phi \equiv \b'\t^{q_t}\s^{p_t} \cdots \t^{q_1}\s^{p_1},
\tag 2.1
$$where $\b'$ is
a Whitehead automorphism of $F_2$ of type (W1), $t \in \Bbb {N}$
and both $p_i, \, q_i$ are (not necessarily positive) integers for
every $i=1, \dots, t$. If not every $p_i$ and $q_i$ has the same
sign (including $0$), apply repeatedly Lemma ~2.2 to the chain on
the right-hand side of (2.1) to obtain that either $\phi \equiv
\b'\pi^r \t^{m_k}\s^{l_k} \cdots \t^{m_1}\s^{l_1}$ or $\phi \equiv
\b'\pi^r \t^{-m_k}\s^{-l_k} \cdots \t^{-m_1}\s^{-l_1}$, where
$\pi$ is as in Lemma ~2.2, $r \in \Bbb Z$, $k \in \Bbb {N}$, and
both $l_i, \, m_i \ge 0$ for every $i=1, \dots, k$. Putting
$\b=\b'\pi^r$, we get the required result. \qed
\enddemo

Under the same notation as in the statement of Lemma ~2.3, we
define the {\it length of an automorphism} $\phi$ of $F_2$ as
$\sum_{i=1}^k (m_i+l_i)$, which is denoted by $|\phi|$. Then
obviously $|\phi|=|\phi'|$.

\proclaim {Lemma 2.4} Let $u, \, v$ be elements in $F_2$. Also let
$m$ be an arbitrary positive integer, and let $\Lambda$ be the set
of all chains of the form (C1) or (C2) of length less than or
equal to $m$. Suppose that $\|\psi(u)\|=\|\psi(v)\|$ for every
$\psi \in \Lambda$. Then we have both $n([\psi(u)];
x)=n([\psi(v)]; x)$ and $n([\psi(u)]; y)=n([\psi(v)]; y)$ for
every $\psi \in \Lambda$.
\endproclaim

\demo{Proof} Under the given hypothesis of the lemma, [2, Lemma
~2.2] yields that $n([u];x)=n([v];x)$ and $n([u];y)=n([v];y)$,
thus proving the lemma when $\psi=1$. Now assuming that the
assertion of the lemma is true for every $\psi_1 \in \Lambda$ with
$|\psi_1| = m' < m$, we shall prove that $n([\psi_2(u)];
x)=n([\psi_2(v)]; x)$ and $n([\psi_2(u)]; y)=n([\psi_2(v)]; y)$
for every $\psi_2 \in \Lambda$ with $|\psi_2|=m'+1$. Such $\psi_2$
can be expressed as $\s^{\pm 1} \psi_1$ or $\t^{\pm 1} \psi_1$ for
some $\psi_1 \in \Lambda$ with $|\psi_1| = m'$.

First let $\psi_2=\s^{\pm 1} \psi_1$. Then clearly
$n([\psi_2(u)];x)=n([\psi_1(u)];x)$ and
$n([\psi_2(v)];x)=n([\psi_1(v)];x)$. Since
$n([\psi_1(u)];x)=n([\psi_1(v)];x)$ by the induction hypothesis,
we have $n([\psi_2(u)]; x)=n([\psi_2(v)]; x)$. Moreover it is
clear that $n([\psi_2(u)]; y)=\|\psi_2(u)\|-n([\psi_2(u)]; x)$ and
$n([\psi_2(v)]; y)=\|\psi_2(v)\|-n([\psi_2(v)]; x)$. Since
$\|\psi_2(u)\|=\|\psi_2(v)\|$ by the hypothesis of the lemma, we
finally have $n([\psi_2(u)]; y)=n([\psi_2(v)]; y)$.

The other case where $\psi_2=\t^{\pm 1} \psi_1$ is similar. \qed
\enddemo

For a cyclic word $w$ in $F_2$ and a Whitehead automorphism, say
$\s$, of $F_2$, a subword of the form $xy^rx^{-1}$ ($r \neq 0$),
if any, in $w$ is invariant in passing from $w$ to $\s(w)$,
although there occurs cancellation in $\s(xy^rx^{-1})$ (note that
$\s(xy^rx^{-1})=xy \cdot y^r \cdot y^{-1}x^{-1}=xy^rx^{-1}$). Such
cancellation is called {\it trivial cancellation}. And
cancellation which is not trivial cancellation is called {\it
proper cancellation}. For example, a subword $xy^{-r}x$ ($r \ge
1$), if any, in $w$ is transformed to $xy^{-r+1}xy$ by applying
$\s$, and the cancellation occurring in $\s(xy^{-r}x)$ is proper
cancellation.

\proclaim{Lemma 2.5}  Let $w$ be a cyclic word in $F_2$, and let
$\psi$ be a chain of the form (C1) (or (C2)). If $\psi$ contains
at least $\|w\|$ factors of $\s$ (or $\s^{-1}$), then there cannot
occur proper cancellation in passing from $\psi(w)$ to $\s\psi(w)$
(or $\psi(w)$ to $\s^{-1}\psi(w)$); if $\psi$ contains at least
$\|w\|$ factors of $\t$ (or $\t^{-1}$), then there cannot occur
proper cancellation in passing from $\psi(w)$ to $\t\psi(w)$ (or
$\psi(w)$ to $\t^{-1}\psi(w)$).
\endproclaim

\demo{Proof} We shall show that if $\psi$ is a chain of the form
(C1) such that $\psi$ contains at least $\|w\|$ factors of $\s$,
then no proper cancellation occurs in passing from $\psi(w)$ to
$\s\psi(w)$ (the other cases are similar). Supposing that there is
a chain $\psi'$ of type (C1) such that no proper cancellation
occurs in passing from $\psi'(w)$ to $\s \psi'(w)$, we see that
proper cancellation cannot occur in passing from $\s \psi'(w)$ to
$\s^2 \psi'(w)$ or in passing from $\t^t \s \psi'(w)$ to $\s\t^t\s
\psi'(w)$ for any $t \ge 1$. Hence if there was proper
cancellation in passing from $\psi(w)$ to $\s \psi(w)$, then
proper cancellation would also occur at every step of applying
$\s$ in $\psi$. However since cancelled $y^{\pm 1}$ in proper
cancellation at every step of applying $\s$ in the chain $\s\psi$
must originally exist in $w$ and since the chain $\s\psi$ contains
more than $\|w\|$ factors of $\s$, we reach a contradiction. \qed
\enddemo

\heading 3. Proof of Theorem 1.2
\endheading

We shall prove the following.

\roster \item"($\ast$)" Let $\O$ be the set of all chains of the
form (C1) or (C2) of length less than or equal to $2\|u\|+3$.
Suppose that $\|\psi(u)\|=\|\psi(v)\|$ for every $\psi \in \O$.
Then $u$ and $v$ are translation equivalent in $F_2$.
\endroster

Once ($\ast$) is proved, the translation equivalence
of $u, \, v$ in $F_2$ is algorithmically decidable as follows.

\proclaim {Algorithm} Let $\O$ be the set of all chains of the
form (C1) or (C2) of length less than or equal to $2\|u\|+3$ (note
that $\O$ is a finite set). Check if $\|\psi(u)\|=\|\psi(v)\|$ for
every $\psi \in \O$. If so, conclude that $u$ and $v$ are
translation equivalent in $F_2$; otherwise conclude that $u$ and
$v$ are not translation equivalent in $F_2$.
\endproclaim

Let $\phi$ be an automorphism of $F_2$. By Lemma ~2.3, $\phi$ can
be represented as $\phi \equiv \b\phi'$, where $\b$ is a Whitehead
automorphism of $F_2$ of type (W1) and $\phi'$ is of the form
either (C1) or (C2). We proceed with the proof of ($\ast$) by
induction on $|\phi'|$. Suppose that $\phi'$ is a chain of the
form (C1) with $|\phi'|>2\|u\|+3$ (the case for (C2) is similar).
Assuming that $\|\psi (u)\|=\|\psi(v)\|$ for every chain $\psi$ of
the form (C1) or (C2) with $|\psi| < |\phi'|$, we shall show that
$\|\phi' (u)\|=\|\phi'(v)\|$, which is equivalent to showing that
$\|\phi (u)\|=\|\phi(v)\|$. Suppose that $\phi'$ ends with $\t$
(the case where $\phi'$ ends with $\s$ is similar), that is,
$$\phi'=\t^{m_k}\s^{l_k} \cdots \t^{m_1}\s^{l_1},$$
where both $l_i, \, m_i \ge 0$ for every $i=1, \dots, k$ and
$m_k>0$. Put
$$\phi_1=\t^{m_k-1}\s^{l_k} \cdots \t^{m_1}\s^{l_1}.$$ Also put
$$u_1=\phi_1(u) \quad \text{and} \quad v_1=\phi_1(v).$$ It then follows
from $\t(u_1)=\phi'(u)$ and $\t(v_1)=\phi'(v)$ that
$$\aligned
\|\phi'(u)\|&=\|u_1\|+n([u_1]; y)-2n([u_1]; y,x^{-1}) \\
\|\phi'(v)\|&=\|v_1\|+n([v_1]; y)-2n([v_1]; y,x^{-1}).
\endaligned
\tag 3.1
$$
By the induction hypothesis, we have $\|u_1\|=\|v_1\|$. Moreover,
by Lemma ~2.4, we have $n([u_1]; y)=n([v_1]; y)$. So it suffices
to show $n([u_1]; y,x^{-1})=n([v_1]; y,x^{-1})$ to get the
equality $\|\phi' (u)\|=\|\phi'(v)\|$.

Clearly the chain $\phi_1$ has length $|\phi_1| = |\phi'|-1 \ge
2\|u\|+3$. Hence either $\s$ or $\t$ occurs at least $\|u\|+2$
times in $\phi_1$. We consider two cases accordingly.

\proclaim {Case 1} $\s$ occurs at least $\|u\|+2$ times in
$\phi_1$.
\endproclaim

In this case, clearly $l_k>0$. Put $$u_2=\t^{m_k-1}\s^{l_k-1}
\cdots \t^{m_1}\s^{l_1}(u) \quad \text{and} \quad u_2'=\s^{l_k-1}
\cdots \t^{m_1}\s^{l_1}(u).$$ Then $u_2=\t^{m_k-1}(u_2')$. In the
following claims, we shall make some observations about the cyclic
word $[u_2']$.

\proclaim {Claim 1} (i) If $l_k-1>0$, then $[u_2']$ does not have
$x^2$ or $x^{-2}$ as a subword.

(ii) Let $l_k-1=0$. Then the cyclic word $[\s^{l_{(k-1)}} \cdots
\t^{m_1}\s^{l_1}(u)]$ does not have $x^2$ or $x^{-2}$ as a
subword. If there is a subword $x^2$ or $x^{-2}$ in $[u_2']$, then
it is actually part of the subword $yx^2$ or $x^{-2}y^{-1}$,
respectively.
\endproclaim

\demo{Proof of Claim 1} (i) Let $l_k-1>0$. Since the chain
$\s^{l_k-2} \cdots \t^{m_1}\s^{l_1}$ contains at least $\|u\|$
factors of $\s$, by Lemma ~2.5 no proper cancellation occurs in
passing from $[\s^{l_k-2} \cdots \t^{m_1}\s^{l_1}(u)]$ to
$[\s^{l_k-1} \cdots \t^{m_1}\s^{l_1}(u)]=[u_2']$. This yields that
$x^2$ or $x^{-2}$ cannot occur in $[u_2']$ as a subword.

(ii) Let $l_k-1=0$. Then $l_{(k-1)}>0$ and the chain
$\s^{l_{(k-1)}-1} \cdots \t^{m_1}\s^{l_1}$ contains at least
$\|u\|$ factors of $\s$. Again by Lemma ~2.5, no proper
cancellation occurs in passing from $[\s^{l_{(k-1)}-1} \cdots
\t^{m_1}\s^{l_1}(u)]$ to $[\s^{l_{(k-1)}} \cdots
\t^{m_1}\s^{l_1}(u)]$. This yields that $x^2$ or $x^{-2}$ cannot
occur in $[\s^{l_{(k-1)}} \cdots \t^{m_1}\s^{l_1}(u)]$ as a
subword.

Thus if there exists $x^2$ or $x^{-2}$ in $[u_2']$ as a subword,
it must have newly occurred in passing from $[\s^{l_{(k-1)}}
\cdots \t^{m_1}\s^{l_1}(u)]$ to $[\t^{m_{(k-1)}}\s^{l_{(k-1)}}
\cdots \t^{m_1}\s^{l_1}(u)]=[u_2']$. This implies that if there is
a subword $x^2$ or $x^{-2}$ in $[u_2']$, it is actually part of
the subword $yx^2$ or $x^{-2}y^{-1}$, respectively. \qed
\enddemo

\proclaim {Claim 2} The cyclic word $[u_2']$ can be written as
$[w_1z_1 \cdots w_tz_t]$, where $z_i$ is either $xy^tx^{-1}$ or
$xy^{-t}x^{-1}$ ($t \ge 1$), and $w_i$ contains no $yx^{-1}$ or
$xy^{-1}$ as a subword and neither begins with nor ends with
$x^{\pm 1}$.
\endproclaim

\demo{Proof of Claim 2} Since the chain $\s^{l_k-1} \cdots
\t^{m_1}\s^{l_1}$ contains at least $\|u\|+1$ factors of $\s$, by
Lemma ~2.5 no proper cancellation occurs in passing from $[u_2']$
to $[\s(u_2')]$. This implies that any subword $yx^{-1}$ or
$xy^{-1}$, if any, in $[u_2']$ must be part of a subword of the
form $xy^tx^{-1}$ or $xy^{-t}x^{-1}$ ($t \ge 1$), respectively, in
$[u_2']$.

Suppose that $xy^tx^{-2}$ or $x^2y^{-t}x^{-1}$ ($t \ge 1$) occurs
in $[u_2']$ as a subword. By Claim ~1 ~(i), this happens only when
$l_k-1=0$. Also by the second part of Claim ~1 ~(ii), any subword
of the form $xy^tx^{-2}$ or $x^2y^{-t}x^{-1}$ ($t \ge 1$) in
$[u_2']$ is part of a subword of the form $xy^tx^{-s}y^{-1}$ or
$yx^sy^{-t}x^{-1}$ ($s \ge 2$), respectively, in $[u_2']$. But
then a subword of the form $yx^{-s}y^{-1}$ or $yx^sy^{-1}$ ($s \ge
2$) must exist in $[\s^{l_{(k-1)}} \cdots \t^{m_1}\s^{l_1}(u)]$, a
contradiction to the first part of Claim ~1 ~(ii). \qed
\enddemo

Now put $$u_1'=\s(u_2').$$ By Claim ~2, we have $[u_1']=[\s(w_1z_1
\cdots w_tz_t)]= [w_1'z_1 \cdots w_t'z_t]$, where $w_i'=(\{x\},
y)(w_i)$. Then $w_i'$ contains no $yx^{-1}$ or $xy^{-1}$ as a
subword and has the same initial and terminal letters as $w_i$
does for each $i$. Since $u_1$ and $u_2$ are obtained by applying
$\t^{m_k-1}$ to $u_1'$ and $u_2'$, respectively, we see that
$$n([u_1]; y,x^{-1})=n([u_2]; y,x^{-1}).$$
Arguing similarly, we have
$$n([v_1]; y,x^{-1})=n([v_2]; y,x^{-1}),$$
where $v_2=\t^{m_k-1}\s^{l_k-1} \cdots \t^{m_1}\s^{l_1}(v).$
Furthermore, since
$$\aligned
-2n([u_2]; y,x^{-1})&=\|\t^{m_k}\s^{l_k-1} \cdots
\t^{m_1}\s^{l_1}(u)\| -\|u_2\|-n([u_2];y)\\
-2n([v_2]; y,x^{-1})&=\|\t^{m_k}\s^{l_k-1} \cdots
\t^{m_1}\s^{l_1}(v)\| -\|v_2\|-n([v_2];y),
\endaligned
$$
by the induction hypothesis applied to both $\|\t^{m_k}\s^{l_k-1}
\cdots \t^{m_1}\s^{l_1}(u)\|= \|\t^{m_k}\s^{l_k-1} \cdots
\t^{m_1}\s^{l_1}(v)\|$ and $\|u_2\|=\|v_2\|$ together with Lemma
~2.4 applied to $n([u_2];y)=n([v_2];y)$, we finally have
$$n([u_1]; y,x^{-1})=n([u_2]; y,x^{-1})=n([v_2]; y,x^{-1})=n([v_1]; y,x^{-1}),$$
that is, $n([u_1]; y,x^{-1})=n([v_1]; y,x^{-1}),$ as required.

\proclaim {Case 2} $\t$ occurs at least $\|u\|+2$ times in
$\phi_1$.
\endproclaim

We divide this case into two subcases.

\proclaim {Case 2.1} $m_k \ge 2$.
\endproclaim

Put $$u_3=\t^{m_k-2}\s^{l_k} \cdots \t^{m_1}\s^{l_1}(u) \quad
\text{and} \quad v_3=\t^{m_k-2}\s^{l_k} \cdots
\t^{m_1}\s^{l_1}(v).$$ Here since the chain $\t^{m_k-2}\s^{l_k}
\cdots \t^{m_1}\s^{l_1}$ contains at least $\|u\|+1$ factors of
$\t$, by Lemma ~2.5 no proper cancellation occurs in passing from
$[u_3]$ to $[\t(u_3)]=[u_1]$. Hence we have $n([u_1];
y,x^{-1})=n([u_3]; y,x^{-1})$. Similarly $n([v_1];
y,x^{-1})=n([v_3]; y,x^{-1})$. Since
$$
\aligned -2n([u_3]; y,x^{-1})&=\|\t^{m_k-1}\s^{l_k} \cdots
\t^{m_1}\s^{l_1}(u)\| -\|u_3\|-n([u_3];y) \\
-2n([v_3]; y,x^{-1})&=\|\t^{m_k-1}\s^{l_k} \cdots
\t^{m_1}\s^{l_1}(u)\| -\|v_3\|-n([v_3];y),
\endaligned
$$
the desired equality $n([u_1]; y,x^{-1})=n([v_1]; y,x^{-1})$
follows from the induction hypothesis and Lemma ~2.4.

\proclaim {Case 2.2} $m_k =1$.
\endproclaim

In this case clearly $m_{(k-1)}>0$. Put $$u_4=\s^{l_k}
\t^{m_{(k-1)}-1} \cdots \t^{m_1}\s^{l_1}(u) \quad \text{and} \quad
v_4=\s^{l_k} \t^{m_{(k-1)}-1} \cdots \t^{m_1}\s^{l_1}(v).$$ As in
Case 1, we can see that
$$\aligned
n([u_1]; y,x^{-1})&=n([u_4]; y,x^{-1})\\
n([v_1]; y,x^{-1})&=n([v_4]; y,x^{-1}).
\endaligned
$$
Then since
$$\aligned
-2n([u_4]; y,x^{-1})&=\|\t \s^{l_k} \t^{m_{(k-1)}-1} \cdots
\t^{m_1}\s^{l_1}(u)\| -\|u_4\|-n([u_4];y) \\
-2n([v_4]; y,x^{-1})&=\|\t \s^{l_k} \t^{m_{(k-1)}-1} \cdots
\t^{m_1}\s^{l_1}(u)\| -\|v_4\|-n([v_4];y),
\endaligned
$$
the required equality $n([u_1]; y,x^{-1})=n([v_1]; y,x^{-1})$
follows from the induction hypothesis and Lemma ~2.4.

The proof of ($\ast$), and hence the proof of Theorem ~1.2, is now
completed. \qed

\heading Acknowledgements
\endheading
The author is grateful to Ilya Kapovich and Vladimir Shpilrain for
suggesting this research topic. The author is also thankful to the
referee for a careful report. This work was supported by the Korea
Research Foundation Grant funded by the Korean Government
(KRF-2006-531-C00011).

\heading References
\endheading

\roster

\item"1." G. Baumslag, A. G. Myasnikov and V. Shpilrain, {\it Open
problems in combinatorial group theory, Second edition}, Contemp.
Math. {\bf 296} (2002), 1--38.

\item"2." I. Kapovich, G. Levitt, P. E. Schupp and V. Shpilrain,
Translation equivalence in free groups, {\it Trans. Amer. Math.
Soc.}, to appear.

\item"3." D. Lee, Translation equivalent elements in free groups,
{\it J. Group Theory}, to appear.

\item"4." D. Lee, Counting words of minimum length in an
automorphic orbit, {\it J. Algebra} {\bf 301} (2006), 35--58.

\item"5." C. J. Leininger, Equivalent curves in surfaces, {\it
Geom. Dedicata} {\bf 102} (2003), 151--177.

\item"6." J. H. C. Whitehead, Equivalent sets of elements in a
free group, {\it Ann. of Math.} {\bf 37} (1936), 782--800.
\endroster
\enddocument